\documentclass[11pt,a4paper]{amsart}
\usepackage{setspace}
\usepackage{amssymb}
\usepackage[all]{xy}
\newtheorem{thm1}{Theorem}[section]

\newtheorem{rem1}[thm1]{Remark}
\newtheorem{def1}[thm1]{Definition}
\newtheorem{cor1}[thm1]{Corollary}

\newtheorem{prop1}[thm1]{Proposition}

\begin{document}
\title[Arithmetical rank of lattice ideals]
{Specializations of multigradings and the arithmetical rank of lattice ideals}
\author[A. Katsabekis]{Anargyros Katsabekis}
\address{Department of Mathematics, University of the Aegean,
83200 Karlovassi, Samos, GREECE} \email{katsabek@aegean.gr}
\author[A. Thoma]{ Apostolos Thoma}
\address{Department of Mathematics, University of Ioannina,
Ioannina 45110, GREECE } \email{athoma@uoi.gr}


\begin{abstract}

\par In this article we study specializations of multigradings
and apply them to the problem of the computation of the
arithmetical rank of a lattice ideal $I_{L_{\mathcal{G}}} \subset
K[x_{1},\ldots,x_{n}]$. The
arithmetical rank of $I_{L_{\mathcal{G}}}$ equals the
$\mathcal{F}$-homogeneous arithmetical rank of
$I_{L_{\mathcal{G}}}$, for an appropriate specialization
$\mathcal{F}$ of $\mathcal{G}$. To the lattice ideal
$I_{L_{\mathcal{G}}}$ and every specialization $\mathcal{F}$ of
$\mathcal{G}$ we associate a simplicial complex. We prove that
combinatorial invariants of the simplicial complex provide lower
bounds for
 the $\mathcal{F}$-homogeneous
arithmetical rank of $I_{L_{\mathcal{G}}}$.
\end{abstract}
\maketitle

\section{Introduction }

\par On the polynomial ring  $S=K[x_1,\dots, x_n]$ with coefficients in a field $K$ one can impose
several multigradings defined by abelian groups. Let $\mathcal{G}$
be a finitely generated abelian group together with a
distinguished ordered set $\{{\bf g}_1,\ldots,{\bf g}_n\}$ of $n$
generators. The degree map
$${\rm deg}_{\mathcal{G}}: \mathbb{Z}^{n} \rightarrow \mathcal{G}, \
{\rm deg}_{\mathcal{G}}({\bf u})=u_{1}{\bf g}_1+\cdots+u_{n}{\bf
g}_n \ \textrm{for} \ {\bf u}=(u_{1},\ldots,u_{n}) \in
\mathbb{Z}^n,$$ defines a multigrading on $S$ by $\mathcal{G}$.
The $\mathcal{G}$-{\em degree} of the monomial $x_1^{u_1} \cdots
x_n^{u_n}$ is ${\rm deg}_{\mathcal{G}}({\bf u})$. A polynomial $F
\in S$ is called $\mathcal{G}$-{\em homogeneous} if the monomials
in each non zero term of $F$ have the same $\mathcal{G}$-degree.
An ideal $J$ is called $\mathcal{G}$-homogeneous if it is
generated by $\mathcal{G}$-homogeneous polynomials.\\The grading
on $S$ by $\mathcal{G}$ defines the exact sequence $$0
\longrightarrow L_{\mathcal{G}} \stackrel{i}{\longrightarrow}
\mathbb{Z}^{n} \stackrel{{\rm deg}_{\mathcal{G}}}{\longrightarrow}
\mathcal{G} \longrightarrow 0.$$ Depending on the emphasis given
to the group $\mathcal{G}$ or the lattice $L_{\mathcal{G}}$ it is
called $\mathcal{G}$-grading or $L_{\mathcal{G}}$-grading. Remark
that $\mathcal{G}$ together with the set $\{{\bf g}_1,\ldots,{\bf
g}_n\}$ determines the lattice $$L_{\mathcal{G}}=\{{\bf u} \in
\mathbb{Z}^n: {\rm deg}_{\mathcal{G}}({\bf
u})=\mathbf{0}_{\mathcal{G}}\} \subset \mathbb{Z}^n$$ of relations
of ${\bf g}_1,\ldots,{\bf g}_n$. A lattice $L \subset
\mathbb{Z}^n$ determines the group $\mathcal{G}=\mathbb{Z}^{n}/L$
and a distinguished set of $n$ generators ${\bf g}_i={\bf e}_i+L$
for every $i=1,\ldots,n$, where ${\bf e}_1,\ldots,{\bf e}_n$ are
the unit vectors of $\mathbb{Z}^{n}$.\\ Multigradings of polynomial
rings have been extensively studied and systematically used over
the last years, see \cite{KR} chapter 4, \cite{MS} chapter 8,
\cite{St} chapter 10. Several times one has to consider coarser
gradings for an $S$-module than the finest one, see \cite{EMS},
\cite{K}, \cite{KM}. This procedure of passing from a finer to a
coarser grading is called specialization or coarsening the
grading, see \cite{KM}. This is the case studied in the present paper.
We are interested in the problem of computing the arithmetical
rank of a toric or lattice ideal. The {\em arithmetical rank},
denoted by ${\rm ara}(J)$, of an ideal $J \subset
K[x_{1},\ldots,x_{n}]$ is the smallest integer $s$ for which there
exist polynomials $F_1,\dots ,F_s$ in $J$ such that $rad(J) =
rad(F_1,\ldots ,F_s)$. The computation of the arithmetical rank of
a lattice ideal $I_{L}$ is a difficult problem and remains open
even in very simple cases like the ideal of the Macaulay curve in
the three dimensional projective space, see \cite{E} chapter 15.
Every lattice ideal $I_{L}$ has a natural multigraded structure
(\cite{C}, \cite{St}), in fact it is $\mathcal{G}$-homogeneous for
$\mathcal{G}=\mathbb{Z}^{n}/L$. The lattice ideal $rad(I_{L})$ can always be generated up to radical by $\mathcal{G}$-homogeneous
polynomials, and sometimes this is possible with ${\rm ara}(I_{L})$ such polynomials, as was shown in \cite{Eto}, \cite{K}. But this is not the case
in general. In an example of a lattice ideal studied in \cite{KMT}
the arithmetical rank is somewhere between $80$ to $90$ while the
minimum number of $\mathcal{G}$-homogeneous polynomials needed to
generate $rad(I_{L})$ up to radical is exactly $1740$. This means that
$\mathcal{G}$-homogeneous polynomials are not always enough to
minimally generate the radical of a lattice ideal up to radical.
Therefore one has to better understand non
$\mathcal{G}$-homogeneous set-theoretic intersections for lattice
ideals. A first step in this direction is to consider coarser
$\mathcal{F}$-gradings than the $\mathcal{G}$-grading and study
the minimum number of $\mathcal{F}$-homogeneous polynomials needed
to generate the radical of a lattice ideal up to radical. We will define a
relation $\leq$ on the set of gradings by groups with $n$
generators. The grading defined by a group $\mathcal{F}=<{\bf
f}_1,\ldots,{\bf f}_n>$ is called a {\em specialization of}
$\mathcal{G}=<{\bf g}_1,\ldots,{\bf g}_n>$ if every
$\mathcal{G}$-homogeneous ideal in $K[x_1,\ldots,x_n]$ is also
$\mathcal{F}$-homogeneous and this will be denoted by $\mathcal{F} \leq
\mathcal{G}$. Specializations of $\mathcal{G}$-gradings were used
in \cite{K} to compute concrete polynomial equations that
set-theoretically define certain toric varieties. Section 2 of the
paper is devoted to a more systematic study of specializations of
$\mathcal{G}$-gradings.

Let $\mathcal{F}$ be a specialization of $\mathcal{G}$. Given a
$\mathcal{G}$-homogeneous ideal $J \subset K[x_{1},\ldots,x_n]$,
it is natural to define the $\mathcal{F}$-{\em homogeneous
arithmetical rank} of $J$, denoted by ${\rm
ara}_{\mathcal{F}}(J)$, as the smallest integer $s$ such that
$rad(J) = rad(F_1,\ldots ,F_s)$ and all the polynomials
$F_1,\ldots ,F_s$ are $\mathcal{F}$-homogeneous. This notion is
important for two reasons: \begin{enumerate} \item It is an upper
bound for ${\rm ara}(J)$. More precisely for a
$\mathcal{G}$-homogeneous ideal $J$ and a specialization
$\mathcal{F}$ of $\mathcal{G}$ it holds:
$${\rm ht}(J)\leq {\rm ara}(J) \leq {\rm ara}_{\mathcal{F}}(J)\leq {\rm ara}_{\mathcal{G}}
(J),$$ where ${\rm ht}(J)$ is the height of $J$. When ${\rm
ht}(J)={\rm ara}(J)$ the ideal $J$ is called {\em set-theoretic
complete intersection} and when ${\rm ht}(J)={\rm
ara}_{\mathcal{F}}(J)$ it is called $\mathcal{F}$-{\em homogeneous
set-theoretic complete intersection}. \item For every
$\mathcal{G}$-homogeneous ideal $J$ there is an
$\mathcal{F}$-grading such that ${\rm ara}(J)={\rm
ara}_{\mathcal{F}}(J)$ and $\mathcal{F} \leq \mathcal{G}$, see
Proposition 3.3. \end{enumerate} The most difficult part in
computing the arithmetical rank or $\mathcal{F}$-homogeneous
arithmetical rank of a lattice ideal $I_{L_{\mathcal{G}}}$ is to
find sharp lower bounds. Lower bounds of the arithmetical rank of
$I_{L_{\mathcal{G}}}$ can be provided sometimes by local or
etal{\'e} cohomology, see \cite{B}, \cite{BL}. The main result of
this article, Theorem 3.10, generalizes the results of \cite{KMT},
\cite{KT} and provides lower bounds for the
$\mathcal{F}$-homogeneous arithmetical rank of the lattice ideal
$I_{L_{\mathcal{G}}}$, where $\mathcal{F}$ is a specialization of
$\mathcal{G}$, using combinatorial invariants of a simplicial
complex associated to the ideal $I_{L_{\mathcal{G}}}$ and the
specialization $\mathcal{F}$ of $\mathcal{G}$. As an application in Section 4 we study an example of a lattice ideal
$I_{L_{\mathcal{G}}}$. We compute the bounds given in Theorem 3.10
and prove that they are sharp. Finally we show that the lattice
ideal $I_{L_{\mathcal{G}}}$ is not a $\mathcal{F}$-homogeneous
set-theoretic complete intersection for infinitely many
specializations $\mathcal{F}$ of $\mathcal{G}$.\\
\newline

\section{Basic theory of specializations of $\mathcal{G}$-gradings}

\subsection{Preliminaries}
\par Given a
lattice $L \subset \mathbb{Z}^n$, the ideal
$$I_{L}:=(\{{\bf x}^{{\bf \alpha }_+}-{\bf x}^{{\bf \alpha } _-} \ | \ {\bf \alpha } ={\bf \alpha } _+-{\bf
\alpha} _- \in L\})\subset K[x_1,\dots ,x_n]$$ is called {\em
lattice ideal}. Here ${\bf \alpha} _+\in \mathbb{N}^n$ and ${\bf
\alpha } _-\in \mathbb{N}^n$ denote the positive and negative part of
${\bf \alpha } $, respectively, and ${\bf x}^{{\bf \beta
}}=x_1^{b_1}\cdots x_n^{b_n}$ for ${\bf \beta } =(b_1,\dots
,b_n)\in \mathbb{N}^n$. The {\em saturation} of a sublattice $L$ of
$\mathbb{Z}^n$ is the lattice
$$ Sat(L):=\{ {\bf \alpha } \in \mathbb{Z}^n \ | \ d {\bf \alpha } \in L
\ \textrm{for some} \ d \in \mathbb{Z}^*\}.$$ We say that the lattice
$L$ is {\em saturated} if $L=Sat(L)$. This is equivalent to saying
that the group $\mathbb{Z}^{n}/L$ is torsion free. The lattice
ideal $I_{L}$ is prime if and only if $L$ is saturated. A prime
lattice ideal is called a {\em toric ideal}, while the set of
zeroes in $K^n$ is an {\em affine toric variety} in the sense of
\cite{St}.\\ If $L=<{\bf l}_1,\ldots,{\bf l}_k>$ is a sublattice
of $\mathbb{Z}^n$ of rank $k<n$, then there exists a set of vectors
$A=\{{\bf a}_1,\ldots,{\bf a}_n\} \subset \mathbb{Z}^m$ such that
$Sat(L)=L_{\mathbb{Z}A}$, where $m=n-k$ and
$\mathbb{Z}A=\{q_{1}{\bf a}_1+ \cdots + q_{n}{\bf a}_n:
q_{1},\ldots,q_{n} \in \mathbb{Z}\}$ is the lattice spanned by
$A$. Remark that $L_{\mathbb{Z}A}$ is saturated. In order to determine $A$ we work as follows. Set ${\bf
L}=({\bf l}_1,\ldots,{\bf l}_k)$ the matrix with columns ${\bf
l}_1,\ldots,{\bf l}_k$, then there are unimodular integral
matrices ${\bf U}$ and ${\bf Q}$ of orders $n$ and $k$,
respectively, such that ${\bf U}{\bf L}{\bf
Q}=diag({\lambda}_1,\ldots, {\lambda}_k,0,\ldots,0)$ is in Smith
normal form. Here ${\lambda}_1,\ldots,{\lambda}_k$ are natural
numbers and ${\lambda}_i$ divides ${\lambda}_{i+1}$. The set $A$
can be chosen as the one consisting of the columns of the matrix formed by the last
$n-k$ rows of ${\bf U}$. Moreover the group $\mathbb{Z}^n/L$ is
isomorphic to $\mathbb{Z}^m\oplus \mathbb{Z}_{\lambda_1} \oplus
\cdots \oplus \mathbb{Z}_{\lambda_k}$, \cite{MS}. We can associate
with the lattice ideal $I_{L}$ the rational polyhedral cone
$$\sigma_{A}:= pos_{\mathbb{Q}}(A)=\{\sum_{i=1}^{n}d_{i}{\bf a}_i \
| \ d_{i} \in \mathbb{Q}_{\geq 0}\}.$$ A {\em face} of
$\sigma_{A}$ is any set of the form $$\mathcal{T}=\sigma_{A} \cap
\{{\bf x} \in \mathbb{Q}^m: {\bf c}{\bf x}=0\}$$ where ${\bf c}
\in \mathbb{Q}^{m}$ and ${\bf c}{\bf x} \geq 0$ for all ${\bf x} \in
\sigma_{A}$. Faces of dimension one are called {\em extreme rays}.
A cone $\sigma_{A}$ is {\em strongly convex} if $\{{\bf 0}\}$ is a
face of $\sigma_{A}$, where ${\bf 0}=(0,\ldots,0)$.

\subsection{Specializations of $\mathcal{G}$-gradings}
\par The next theorem indicates that the specialization
property reflects on the lattice of relations of the generators of
$\mathcal{G}$ and correspondingly in the lattice ideal
$I_{L_{\mathcal{G}}}$.

\begin{thm1} Let $\mathcal{F}=<{\bf f}_{1},\ldots,{\bf f}_{n}>$
and $\mathcal{G}=<{\bf g}_{1},\ldots,{\bf g}_{n}>$ be finitely
generated abelian groups. The following are
equivalent:\\
(a) $\mathcal{F} \leq \mathcal{G}$, i.e. $\mathcal{F}$ is a specialization of $\mathcal{G}$.\\
(b) $L_{\mathcal{G}} \subset L_{\mathcal{F}}$.\\ (c)
$I_{L_{\mathcal{G}}} \subset I_{L_{\mathcal{F}}}$.\\ (d) There is
a group epimorphism from $\mathcal{G}$ to $\mathcal{F}$, sending
${\bf g}_i$ to ${\bf f}_i$.
\end{thm1}
\textbf{Proof.} The equivalence (b) $\Leftrightarrow$ (c) is easily derived from the fact that a binomial ${\bf x}^{\bf u}-{\bf x}^{\bf v}$ belongs to a lattice ideal $I_{L_{\mathcal{G}}}$ if and only if the vector ${\bf u}-{\bf v}$ belongs to $L_{\mathcal{G}}$. We will prove that (a) $\Leftrightarrow$ (b).\\
(a) $\Rightarrow$ (b) Let ${\bf u}={\bf u}_{+}-{\bf u}_{-} \in L_{\mathcal{G}}$, where ${\bf u}_{+}=(u_{+,1},\ldots,u_{+,n})$ and ${\bf u}_{-}=(u_{-,1},\ldots,u_{-,n})$. Then $$u_{+,1}{\bf g}_{1}+\cdots+u_{+,n}{\bf g}_{n}=u_{-,1}{\bf g}_{1}+\cdots+u_{-,n}{\bf g}_{n}.$$The ideal $J=({\bf x}^{{\bf u}_{+}}-{\bf x}^{{\bf u}_{-}})$ is $\mathcal{G}$-homogeneous, so it is also $\mathcal{F}$-homogeneous. Thus $$u_{+,1}{\bf f}_{1}+\cdots+u_{+,n}{\bf f}_{n}=u_{-,1}{\bf f}_{1}+\cdots+u_{-,n}{\bf f}_{n}$$ and therefore ${\bf u} \in L_{\mathcal{F}}$.\\
(b) $\Rightarrow$ (a) Let $J \subset K[x_{1},\ldots,x_{n}]$ be a $\mathcal{G}$-homogeneous ideal and ${\bf x}^{\bf u}$, ${\bf x}^{\bf v}$ two monomials of a $\mathcal{G}$-homogeneous generator $F$ of $J$, where ${\bf u}=(u_{1},\ldots,u_{n})$ and ${\bf v}=(v_{1},\ldots,v_{n})$. We have$$u_{1}{\bf g}_{1}+\cdots+u_{n}{\bf g}_{n}=v_{1}{\bf g}_{1}+\cdots+v_{n}{\bf g}_{n},$$ which implies that the vector ${\bf w}=(u_{1}-v_{1},\ldots,u_{n}-v_{n})$ belongs to $L_{\mathcal{G}}$. But $L_{\mathcal{G}} \subset L_{\mathcal{F}}$, so ${\bf w}$ belongs to $L_{\mathcal{F}}$ and therefore $$u_{1}{\bf f}_{1}+\cdots+u_{n}{\bf f}_{n}=v_{1}{\bf f}_{1}+\cdots+v_{n}{\bf f}_{n}.$$Thus $J$ is $\mathcal{F}$-homogeneous.\\
Finally we will prove
that (b) $\Leftrightarrow$ (d). Assume first that $L_{\mathcal{G}}
\subset L_{\mathcal{F}}$. We define $\phi: \mathcal{G} \rightarrow
\mathcal{F}$ by setting $$\phi(\alpha_{1}{\bf
g}_1+\cdots+\alpha_{n}{\bf g}_n)=\alpha_{1}{\bf
f}_1+\cdots+\alpha_{n}{\bf f}_n.$$ The map $\phi$ is well defined.
Let ${\bf u} \in \mathcal{G}$ be such that ${\bf u}=\alpha_{1}{\bf
g}_1+\cdots+\alpha_{n}{\bf g}_n$ and ${\bf u}=\beta_{1}{\bf
g}_1+\cdots+\beta_{n}{\bf g}_n$. Then the vector
$(\alpha_{1}-\beta_{1},\ldots,\alpha_{n}-\beta_{n})$ belongs to
$L_{\mathcal{G}}$, which is a subset of $L_{\mathcal{F}}$ and
therefore $\alpha_{1}{\bf f}_1+\cdots+\alpha_{n}{\bf
f}_n=\beta_{1}{\bf f}_1+\cdots+\beta_{n}{\bf f}_n$. Obviously
$\phi$ is a homomorphism mapping $\mathcal{G}$ onto
$\mathcal{F}$.\\Conversely assume that there is a group
epimorphism $\phi: \mathcal{G} \rightarrow \mathcal{F}$, sending
${\bf g}_i$ to ${\bf f}_i$. Let ${\bf u}=(u_1,\ldots,u_n) \in
L_{\mathcal{G}}$, then $u_1{\bf g}_1+\cdots+u_n{\bf
g}_n=\mathbf{0}_{\mathcal{G}}$ and therefore $\phi(u_1{\bf
g}_1+\cdots+u_n{\bf g}_n)=\mathbf{0}_{\mathcal{F}}$. Thus $u_1{\bf
f}_1+\cdots+u_n{\bf f}_n=\mathbf{0}_{\mathcal{F}}$, which implies
that ${\bf u}$ belongs to $L_{\mathcal{F}}$. \hfill $\square$\\

\begin{rem1} {\rm Let $\phi $ be the epimorphism defined in the proof of Theorem 2.1. Any
$\mathcal{G}$-graded $S$-module $M$ can be regarded as an
$\mathcal{F}$-graded module by setting $M_{\bf u}=\bigoplus_{{\bf
v} \in \phi^{-1}({\bf u})} M_{\bf v}$.} \end{rem1}

\begin{cor1} Let $\mathcal{F}$, $\mathcal{G}$ be finitely generated abelian groups with $n$
generators and $A$, $B$ sets of vectors such that
$Sat(L_{\mathcal{F}})=L_{\mathbb{Z}B}$ and
$Sat(L_{\mathcal{G}})=L_{\mathbb{Z}A}$. If $\mathcal{F}$ is a
specialization of $\mathcal{G}$, then $\mathbb{Z}B$ is a
specialization of $\mathbb{Z}A$.
\end{cor1}

\textbf{Proof.} From Theorem 2.1 we have that $L_{\mathcal{G}} \subset
L_{\mathcal{F}}$ and therefore $Sat(L_{\mathcal{G}}) \subset
Sat(L_{\mathcal{F}})$. Thus $L_{\mathbb{Z}A} \subset
L_{\mathbb{Z}B}$, so, again from Theorem 2.1, the group
$\mathbb{Z}B$ is a specialization of $\mathbb{Z}A$. \hfill
$\square$

\begin{rem1} {\rm The group $\mathbb{Z}B$ is a specialization of $\mathcal{F}$, since $L_{\mathcal{F}} \subset Sat(L_{\mathcal{F}})$, and similarly $\mathbb{Z}A$ is a specialization of $\mathcal{G}$.}
\end{rem1}
Let $\widehat{\pi}:\mathbb{Q}^{m} \rightarrow \mathbb{Q}^{r}$ be a
rational affine map with $\widehat{\pi}(\sigma_{A})=\sigma_{B}$.
The restriction
$$\pi:=\widehat{\pi}|_{\sigma_{A}}:\sigma_{A} \rightarrow \sigma_{B}$$ is called {\em projection of cones}.
\begin{prop1} If $\mathbb{Z}B$ is a specialization of $\mathbb{Z}A$, for
$A=\{{\bf a}_1,\ldots,{\bf a}_n\}$ and $B=\{{\bf b}_1,\ldots,{\bf
b}_n\}$, then there is a projection of cones $\pi:\sigma_{A}
\rightarrow \sigma_{B}$ given by $\pi({\bf a}_i)={\bf b}_i$ for
all $i=1,\ldots,n$.
\end{prop1}

\textbf{Proof.} From Theorem 2.1 we have that $L_{\mathbb{Z}A} \subset L_{\mathbb{Z}B}$, since $\mathbb{Z}B$ is a specialization of $\mathbb{Z}A$, so, from Theorem 2.2 in \cite{K}, there is a projection of cones $\pi:\sigma_{A}
\rightarrow \sigma_{B}$ given by $\pi({\bf a}_i)={\bf b}_i$ for
all $i=1,\ldots,n$. \hfill
$\square$\\

We say that $\mathcal{F}$ is {\em equivalent} to $\mathcal{G}$,
denoted by $\mathcal{F} \sim \mathcal{G}$, if every
$\mathcal{F}$-homogeneous ideal is also $\mathcal{G}$-homogeneous
and conversely.
\begin{cor1}  Let $\mathcal{F}=<{\bf f}_{1},\ldots,{\bf f}_{n}>$ and
$\mathcal{G}=<{\bf g}_{1},\ldots,{\bf g}_{n}>$ be finitely
generated abelian groups. The following are
equivalent:\\
(a) $\mathcal{F} \sim \mathcal{G}$.\\
(b) $I_{L_{\mathcal{F}}}=I_{L_{\mathcal{G}}}$.\\ (c)
$L_{\mathcal{F}}= L_{\mathcal{G}}$.\\ (d) $\mathcal{F}$,
$\mathcal{G}$ are isomorphic groups and the isomorphism sends
${\bf g}_i$ to ${\bf f}_i$.
\end{cor1}

Although equivalent gradings defined by different groups provide
exactly the same grading in the polynomial ring, it is interesting
to study them for other reasons, including the fact that they give
different toric sets, see \cite{KT-JA}, which has applications to
Algebraic Statistics, see \cite{GMS}.

From now on $\widetilde{\mathcal{G}}$ will denote the equivalence
class of the group $\mathcal{G}$. By writing
$\widetilde{\mathcal{F}} \leq \widetilde{\mathcal{G}}$ we mean
that for every pair of representatives $\mathcal{F}$ and $\mathcal{G}$ of
$\widetilde{\mathcal{F}}$ and $\widetilde{\mathcal{G}}$,
respectively, it holds $\mathcal{F} \leq \mathcal{G}$. From
Theorem 2.1 it is easily derived that if $\mathcal{F} \leq
\mathcal{G}$ and $\mathcal{G} \leq \mathcal{H}$, then $\mathcal{F}
\leq \mathcal{H}$. So $\leq $ is a partial order on the set of equivalence classes of gradings of groups with $n$ generators with respect to relation $\sim$.
Let $\mathcal{F}$ and $\mathcal{G}$ be groups generated by $n$ elements.
We define the {\em join} of $\widetilde{\mathcal{F}}$ and
$\widetilde{\mathcal{G}}$, denoted by $\widetilde{\mathcal{F} \vee
\mathcal{G}}$, to be the equivalence class of the group
$\mathbb{Z}^n/{(L_{\mathcal{F}} \cap L_{\mathcal{G}})}$. The
{\em meet} of $\widetilde{\mathcal{F}}$ and
$\widetilde{\mathcal{G}}$, denoted by $\widetilde{\mathcal{F}
\wedge \mathcal{G}}$, is defined as the equivalence class of the group
$\mathbb{Z}^n/{(L_{\mathcal{F}}+L_{\mathcal{G}})}$. We have that $\mathcal{F}
\wedge \mathcal{G} \leq \mathcal{F}$ and $\mathcal{F}
\wedge \mathcal{G} \leq \mathcal{G}$, since $L_{\mathcal{F}}+L_{\mathcal{G}}$ contains both $L_{\mathcal{F}}$, $L_{\mathcal{G}}$. Moreover if $\mathcal{H} \leq \mathcal{F}$ and $\mathcal{H} \leq \mathcal{G}$, then $\mathcal{H} \leq \mathcal{F}
\wedge \mathcal{G}$ since $L_{\mathcal{F}}+L_{\mathcal{G}}$ is the smallest sublattice of $\mathbb{Z}^{n}$ containing $L_{\mathcal{F}}$ and $L_{\mathcal{G}}$. The finest grading is given by the abelian group
$\mathbb{Z}^n$ with  generators the vectors ${\bf
e}_i=(0,\ldots,0,1,0,\ldots,0)$, where the $1$ is in the
$i$th position. Note that every finitely generated abelian group
$\mathcal{G}$ is a specialization of $\mathbb{Z}^n$, since
$L_{\mathbb{Z}^{n}}=<{\bf 0}>$. The only
$\mathbb{Z}^n$-homogeneous ideals in $K[x_1,\ldots,x_n]$ are the
monomial ideals, while the coarsest grading is given by the zero
group $\mathcal{O}$ generated by the set of $n$ zero vectors ${\bf
o}_i={\bf 0}$. Note that $\mathcal{O}$ is a specialization of
every abelian group $\mathcal{G}$ with $n$ generators and
$I_{L_{\mathcal{O}}}=<x_1-1,\dots ,x_n-1>$. Every ideal in
$K[x_1,\dots ,x_n]$ is $\mathcal{O}$-homogeneous. So actually
$\mathcal{O} \leq \mathcal{G} \leq \mathbb{Z}^n$.

We say that a $\mathcal{G}$-grading is {\em positive} if
$L_{\mathcal{G}} \cap \mathbb{N}^{n}=\{{\bf 0}\}$. This is
equivalent to saying that the rational polyhedral cone $\sigma_{A}$
is strongly convex. Specializations can be used to give an
equivalent characterization of the positivity condition. For more
equivalent conditions, see \cite{MS}, Chapter 8.

\begin{thm1} Let $\mathcal{G}=<{\bf g}_1,\ldots,{\bf g}_n>$ be a
finitely generated abelian group. The $\mathcal{G}$-grading is
positive if and only if there exists a set $M=\{m_1,\ldots,m_n\}$
of positive integers such that $\mathbb{Z}M$ is a specialization of
$\mathcal{G}$.
\end{thm1}
\textbf{Proof.} Suppose first that there is a set $M=\{m_1,\ldots,m_n\}$
with $\mathbb{Z}M \leq \mathcal{G}$ and that the
$\mathcal{G}$-grading is not positive. Then there is a relation
$$\lambda_{1}{\bf g}_1+\cdots+\lambda_{n}{\bf
g}_n=\mathbf{0}_{\mathcal{G}},$$ where every $\lambda_{i} \in
\mathbb{Z}$ is non negative and there is at least one
$\lambda_{j}$ different from zero. Let $\phi: \mathcal{G}
\rightarrow \mathbb{Z}M$ be the group epimorphism, sending ${\bf
g}_i$ to $m_i$. We have that $\phi(\lambda_{1}{\bf
g}_1+\cdots+\lambda_{n}{\bf g}_n)=0$, so $\lambda_{1}\phi({\bf
g}_1)+\cdots+\lambda_{n}\phi({\bf g}_n)=0$ and therefore
$\lambda_{1}m_{1}+\cdots+\lambda_{n}m_{n}=0$. But
$\lambda_{1}m_1+\cdots+\lambda_{n}m_n>0$, since $m_1,\dots ,m_n$
are positive integers and $\lambda_i$ are non negative with
at least one of them different from zero, a contradiction. Suppose
now that the $\mathcal{G}$-grading is positive, this means that
${\bf 0}$ is a face of the corresponding rational polyhedral cone
$\sigma_{A}$. Thus there is a defining vector ${\bf c}_{\bf 0}$ of
the above face such that ${\bf c}_{\bf 0}{\bf a}_i>0$, for every
$i=1,\ldots,n$. Set $m_i={\bf c}_{\bf 0}{\bf a}_i$, for
$i=1,\ldots,n$, then $M=\mathbb{Z}\{m_1,\ldots,m_n\}$ is
specialization of $\mathbb{Z}A$. Let ${\bf u}=(u_{1},\ldots,u_{n})
\in L_{\mathbb{Z}A}$, then $u_{1}{\bf a}_{1}+ \cdots+u_{n}{\bf
a}_{n}={\bf 0}$. So ${\bf c}_{\bf 0}(u_{1}{\bf a}_{1}+
\cdots+u_{n}{\bf a}_{n})=0$ and therefore $u_{1}({\bf c}_{\bf
0}{\bf a}_1)+ \cdots+u_{n}({\bf c}_{\bf 0}{\bf a}_n)=0$. Thus
$L_{\mathbb{Z}A} \subset L_{\mathbb{Z}M}$. From Remark 2.4 we
deduce that $\mathbb{Z}M$ is a specialization of $\mathcal{G}$.
\hfill $\square$\\

Note that if  $\mathcal{F}$ is a specialization of $\mathcal{G}$
and the $\mathcal{F}$-grading is positive, then the
$\mathcal{G}$-grading is positive.\\
\newline

\section{Arithmetical rank of lattice ideals}

\par In this section the first goal is to prove the existence of an
$\mathcal{H}$-grading such that ${\rm
ara}(I_{L_{\mathcal{G}}})={\rm
ara}_{\mathcal{H}}(I_{L_{\mathcal{G}}})$. After that we will
assign to every pair $(\mathcal{F}, \mathcal{G})$ a simplicial
complex $\mathcal{D}_{\mathcal{F}}^{\mathcal{G}}$ and to every
polynomial $F\in K[x_1,\ldots,x_n]$ a subcomplex of
$\mathcal{D}_{\mathcal{F}}^{\mathcal{G}}$, where $\mathcal{F}$ is
a specialization of $\mathcal{G}$. The second goal is to prove
that if
 $F_1,\ldots,F_s$ are
$\mathcal{F}$-homogeneous polynomials and generate
$rad(I_{L_{\mathcal{G}}})$ up to radical, then each of the
subcomplexes corresponding to the polynomials $F_i$ is a simplex
and their union is a spanning subcomplex of
$\mathcal{D}_{\mathcal{F}}^{\mathcal{G}}$. This will enable us to
provide lower bounds for the $\mathcal{F}$-homogeneous
arithmetical rank based on combinatorial invariants of the
simplicial complex $\mathcal{D}_{\mathcal{F}}^{\mathcal{G}}$.

\begin{thm1} Let $\{F_1,\dots ,F_s\}$ be a set of polynomials in  $K[x_1,\ldots,x_n]$. There
exists a finest $\mathcal{F}$-grading
such that all $F_1,\dots ,F_s$ are $\mathcal{F}$-homogeneous. This grading is unique up to equivalence.
\end{thm1}
\textbf{Proof.} Every polynomial
$F_i\neq 0$ can be written as a finite sum of terms, i.e.
$F_{i}=\sum_{j} c_{ij} {\bf x}^{{\bf u}^{i}_j}$ where $K \ni
c_{ij} \neq 0$. Let $L$ be the lattice generated by all the vectors
${\bf u}_{j}^{i}-{\bf u}_{1}^{i}$, for every $i=1,\ldots,s$. The
polynomials $F_1,\dots ,F_s$ are $\mathcal{F}$-homogeneous for
$\mathcal{F}=\mathbb{Z}^{n}/L$. It remains to prove that
$\mathcal{F}$ is the finest. Suppose that  $F_1,\dots ,F_s$ are
also $\mathcal{G}$-homogeneous, then, for every $i=1,\ldots,s$, we
have that ${\rm deg}_{\mathcal{G}}({\bf u}_{j}^{i})={\rm
deg}_{\mathcal{G}}({\bf u}_{1}^{i})$. So the vectors ${\bf
u}_{j}^{i}-{\bf u}_{1}^{i}$ belong to $L_{\mathcal{G}}$ and
therefore $L_{\mathcal{F}} \subset L_{\mathcal{G}}$. Thus
$\mathcal{G} \leq \mathcal{F}$. This fact also implies that
$\mathcal{F}$ is unique up to equivalence. \hfill $\square$

\begin{cor1} Let $\{F_1,\dots ,F_s\}$ be a set of polynomials in  $K[x_1,\ldots,x_n]$ and let
$\mathcal{G}=<{\bf g}_1,\ldots,{\bf g}_n>$ be a finitely generated
abelian group. There exists a finest
$\mathcal{H}$-grading such that all $F_1,\dots ,F_s$ are
$\mathcal{H}$-homogeneous and also $\mathcal{H} \leq \mathcal{G}$. This grading is unique up to equivalence.
\end{cor1}
\textbf{Proof.} From Theorem 3.1 there
exists a finest $\mathcal{F}$-grading
such that the polynomials $F_1,\dots ,F_s$ are
$\mathcal{F}$-homogeneous. This grading is unique up to equivalence. Let $\mathcal{H}=\mathcal{F} \wedge
\mathcal{G}$ be any representative of the class
$\widetilde{\mathcal{F} \wedge \mathcal{G}}$, then $\mathcal{H}
\leq \mathcal{G}$. Moreover $F_1,\ldots,F_s$ are
$\mathcal{H}$-homogeneous, since the ideal generated by
$F_1,\ldots,F_s$ is $\mathcal{F}$-homogeneous and $\mathcal{H}
\leq \mathcal{F}$. To prove that $\mathcal{H}$ is the finest,
assume that the $F_1,\ldots,F_s$ are $\mathcal{M}$-homogeneous and
also $\mathcal{M} \leq \mathcal{G}$. Then $\mathcal{M} \leq
\mathcal{F}$, so $\mathcal{M} \leq \mathcal{F} \wedge
\mathcal{G}=\mathcal{H}$. \hfill $\square$

\begin{prop1} For any $\mathcal{G}$-homogeneous ideal $J \subset K[x_1,\ldots,x_n]$ there is a finest $\mathcal{H}$-grading such that \begin{enumerate} \item $J$ is $\mathcal{H}$-homogeneous
\item $\mathcal{H} \leq \mathcal{G}$ and \item ${\rm ara}(J)={\rm
ara}_{\mathcal{H}}(J)$.
\end{enumerate}
This grading is unique up to equivalence.
\end{prop1}
\textbf{Proof.} Let ${\rm ara}(J)=s$, which implies that $rad(J)=rad(F_1,\ldots,F_s)$ for some polynomials $F_1,\ldots,F_s$ in $K[x_1,\ldots,x_n]$. From Corollary 3.2 there
exists a finest $\mathcal{H}$-grading
such that $F_1,\ldots,F_s$ are $\mathcal{H}$-homogeneous and
$\mathcal{H} \leq \mathcal{G}$. This grading is unique up to equivalence. It follows that ${\rm ara}(J)={\rm
ara}_{\mathcal{H}}(J)$. \hfill $\square$

The next theorem is an easy consequence of Proposition 3.3, since
every lattice ideal $I_{L_{\mathcal{G}}}$ is
$\mathcal{G}$-homogeneous.
\begin{thm1} For any lattice ideal $I_{L_{\mathcal{G}}} \subset K[x_1,\ldots,x_n]$ there is a unique up to equivalence finest $\mathcal{H}$-grading such that \begin{enumerate} \item $I_{L_{\mathcal{G}}}$ is $\mathcal{H}$-homogeneous
\item $\mathcal{H} \leq \mathcal{G}$ and \item ${\rm
ara}(I_{L_{\mathcal{G}}})={\rm
ara}_{\mathcal{H}}(I_{L_{\mathcal{G}}})$.
\end{enumerate}
\end{thm1}

\par Generally it is difficult to compute a priori the grading $\mathcal{H}$ of Theorem 3.4.
But using the theory of simplicial complexes we can find bounds
for the $\mathcal{F}$-homogeneous arithmetical rank of a lattice
ideal $I_{L}$, in the case where the grading induced by the lattice
$L$ is positive. Also note that in several cases one expects that
the group $\mathcal{H}$ of Theorem 3.4 coincides with $\mathcal{O}$. But even in this case one gets interesting results
from the simplicial complex
$\mathcal{D}_{\mathcal{G}}^{\mathcal{G}}$, see Definition 3.5,
such as a lower bound on the number of monomials in the support of
the polynomials that define the radical, but also to the number of
$\mathcal{F}$-homogeneous components, for various $\mathcal{F}$'s.

Let $\mathcal{G}$ be a finitely generated abelian group with $n$
generators and $\sigma_A$ the rational polyhedral cone associated
with the lattice ideal $I_{L_{\mathcal{G}}} \subset
K[x_1,\ldots,x_n]$, for an appropriate set of vectors $A=\{{\bf
a}_1,\ldots,{\bf a}_n\}$. From now on we shall write
$\sigma_{\mathcal{G}}$ instead of $\sigma_{A}$. The {\em relative
interior} of $\sigma_{\mathcal{G}}$, denoted by
$relint_{\mathbb{Q}}(\sigma_{\mathcal{G}})$,
 is the set of all positive rational linear combinations of ${\bf a}_1,\ldots,{\bf a}_n$. When $\sigma_{\mathcal{G}}$ is strongly convex we have that $\sigma_{\mathcal{G}}
=pos_{\mathbb{Q} }({\bf r}_1, \ldots , {\bf r}_t)$, where $\{{\bf r}_1,
\ldots , {\bf r}_t\}$ is a set of integer vectors, one for each
extreme ray of $\sigma_{\mathcal{G}}$. The vectors ${\bf r}_i$ are
called {\em extreme vectors} of $\sigma_{\mathcal{G}}$. Given a
subset $E$ of $\{1,\dots ,t\}$ we denote by
$\sigma_{\mathcal{G}}(E)$ the subcone $pos_{\mathbb{Q}}({\bf r}_i \ | \
i \ \in E)$ of $\sigma_{\mathcal{G}}$. We are going to deal only
with subcones $\sigma_{\mathcal{G}}(E)$, which are not faces of
the cone $\sigma_{\mathcal{G}}$. They form a poset ordered by
inclusion. Let $\{\sigma_{\mathcal{G}}({\mathbb{E}_1}),\ldots
,{\sigma_{\mathcal{G}}({\mathbb{E}_f}})\}$ be the minimal elements
of this poset, which are called the {\em minimal non faces} of
$\sigma_{\mathcal{G}}$. To every specialization we assign a
simplicial complex $\mathcal{D}_{\mathcal{F}}^{\mathcal{G}}$ that
generalizes the complex
$\Delta_{\sigma}=\mathcal{D}_{\mathcal{G}}^{\mathcal{G}}$ defined
in \cite{KMT} and \cite{KT}.
\begin{def1} {\rm Let $\mathcal{F}$ be a specialization of
$\mathcal{G}$ and $\pi:\sigma_{\mathcal{G}} \rightarrow
\sigma_{\mathcal{F}}$ the corresponding projection of cones. We
define $\mathcal{D}_{\mathcal{F}}^{\mathcal{G}}$ to be the
simplicial complex with vertices $\{\mathbb{E}_1,\dots
,\mathbb{E}_f\}$ such that $T \subset
\{\mathbb{E}_1,\ldots,\mathbb{E}_f\}$ belongs to
$\mathcal{D}_{\mathcal{F}}^{\mathcal{G}}$ if and only if
$$\bigcap_{{\mathbb{E}_i} \in T}relint_{\mathbb{Q}}\left(
\pi(\sigma_{\mathcal{G}}({\mathbb{E}_i}))\right) \neq \emptyset
.$$}
\end{def1}

A subcomplex $H$ of a simplicial complex $\mathcal{D}$ is called a {\em spanning subcomplex} if both have exactly the same set of vertices. The following proposition shows that the simplicial complex
$\mathcal{D}_{\mathcal{G}}^{\mathcal{G}}$ is a spanning subcomplex
of $\mathcal{D}_{\mathcal{F}}^{\mathcal{G}}$.
\begin{prop1} Let $\mathcal{F} \leq \mathcal{G}$ be finitely generated abelian
groups with $n$ generators. Then
$$\mathcal{D}_{\mathcal{G}}^{\mathcal{G}} \subset \mathcal{D}_{\mathcal{F}}^{\mathcal{G}}
\subset \mathcal{D}_{\mathcal{O}}^{\mathcal{G}}$$ where
$\mathcal{O}$ is the group generated by the set of $n$ zero
vectors ${\bf o}_i={\bf 0}$. In fact
$\mathcal{D}_{\mathcal{G}}^{\mathcal{G}}$ is a spanning subcomplex
of $\mathcal{D}_{\mathcal{F}}^{\mathcal{G}}$, the simplicial
complex $\mathcal{D}_{\mathcal{F}}^{\mathcal{G}}$ is a spanning
subcomplex of $\mathcal{D}_{\mathcal{O}}^{\mathcal{G}}$ and
$\mathcal{D}_{\mathcal{O}}^{\mathcal{G}}$ is a simplex.
\end{prop1}

\textbf{Proof.} From the definitions of the three simplicial complexes, all of them have the same set of vertices.
Let $T \in \mathcal{D}_{\mathcal{G}}^{\mathcal{G}}$, then
$$\bigcap_{{\mathbb{E}_i} \in T}relint_{\mathbb{Q}}\left(
\sigma_{\mathcal{G}}({\mathbb{E}_i})\right) \neq \emptyset .$$
Hence there exists a ${\bf x} \in \bigcap_{{\mathbb{E}_i} \in
T}relint_{\mathbb{Q}}\left(
\sigma_{\mathcal{G}}({\mathbb{E}_i})\right)$, which implies that
$\pi({\bf x})$ belongs to $\bigcap_{{\mathbb{E}_i} \in
T}relint_{\mathbb{Q}}\left(
\pi(\sigma_{\mathcal{G}}({\mathbb{E}_i}))\right)$. Consequently
$$\bigcap_{{\mathbb{E}_i} \in T}relint_{\mathbb{Q}}\left(
\pi(\sigma_{\mathcal{G}}({\mathbb{E}_i}))\right) \neq \emptyset
.$$ Thus $T \in \mathcal{D}_{\mathcal{F}}^{\mathcal{G}}$ and
therefore $\mathcal{D}_{\mathcal{G}}^{\mathcal{G}} \subset
\mathcal{D}_{\mathcal{F}}^{\mathcal{G}}$.\\
Let $A=\{{\bf a}_1,\ldots,{\bf a}_n\}$ be a set of vectors such
that $Sat(L_{\mathcal{G}})=L_{\mathbb{Z}A}$ and let $\pi_0$ be the
projection of cones sending ${\bf a}_i$ to ${\bf o}_i$. Then, for
every $j=1,\ldots,f$, we have that
$\pi_0(\sigma_{\mathcal{G}}({\mathbb{E}_j}))$ equals $\{{\bf
0}\}$. So
$$\bigcap_{{\mathbb{E}_i} \in T}relint_{\mathbb{Q}}\left(\{{\bf
0}\}\right) \neq \emptyset, \ \textrm{since} \
relint_{\mathbb{Q}}\left(\{{\bf 0}\}\right)=\{\bf 0\}.$$  Thus
$\mathcal{D}_{\mathcal{O}}^{\mathcal{G}}$ is a simplex, so
$\mathcal{D}_{\mathcal{G}}^{\mathcal{G}}$ and
$\mathcal{D}_{\mathcal{F}}^{\mathcal{G}}$ are subcomplexes of
$\mathcal{D}_{\mathcal{O}}^{\mathcal{G}}$.  \hfill $\square$

To every polynomial in $K[x_1,\ldots ,x_n]$ we are going to assign
a series of simplicial complexes, one for each group $\mathcal{G}$
and a specialization $\mathcal{F}$ of $\mathcal{G}$. Recall that $A=\{{\bf a}_{1},\ldots,{\bf a}_{n}\}$ is a set of vectors such that $Sat(L_{\mathcal{G}})=L_{\mathbb{Z}A}$. Let
$N=x_{i_1}^{n_1}\cdots x_{i_s}^{n_s}$ be a monomial in
$K[x_1,\ldots ,x_n]$. Set $A_N:=\{{\bf a}_{i_1},\ldots, {\bf
a}_{i_s}\}$, the cone of $N$ is
$$cone(N):=\bigcap_{A_N \subset \sigma_{\mathcal{G}}(E)} \sigma_{\mathcal{G}}(E) \subset
\sigma_{\mathcal{G}}.$$ Let $F$ be a polynomial in $K[x_1,\ldots
,x_n]$. We associate with $F$ the induced subcomplex
$\mathcal{D}_{\mathcal{G}}^{\mathcal{G}}(F)$ of
$\mathcal{D}_{\mathcal{G}}^{\mathcal{G}}$ consisting of those
vertices $\mathbb{E}_i$ with the property: there exist a monomial
$N$ in $F$ such that $cone(N)=\sigma_{\mathcal{G}}(\mathbb{E}_i)$.
Let $\mathcal{D}_{\mathcal{F}}^{\mathcal{G}}(F)$ be the subcomplex
of $\mathcal{D}_{\mathcal{F}}^{\mathcal{G}}$ induced on the
vertices of $\mathcal{D}_{\mathcal{G}}^{\mathcal{G}}(F)$.
\begin{thm1} Let $\mathcal{F} \leq \mathcal{G}$ be finitely generated abelian
groups. If $F_1,\ldots,F_s$ generate $rad(I_{L_{\mathcal{G}}})$ up
to radical, then $\bigcup_{i=1}^s
\mathcal{D}_{\mathcal{F}}^{\mathcal{G}}(F_i)$ is a spanning
subcomplex of $\mathcal{D}_{\mathcal{F}}^{\mathcal{G}}$.
\end{thm1}
\textbf{Proof.} Let $\mathbb{E}_i$ be a vertex of $\mathcal{D}_{\mathcal{F}}^{\mathcal{G}}$. Then $\mathbb{E}_i$ is a vertex of $\mathcal{D}_{\mathcal{G}}^{\mathcal{G}}$ and therefore, from Theorem 5.1 in \cite{KMT}, there exists a monomial $N$ in some $F_j$ such that $cone(N)=\sigma_{\mathcal{G}}(\mathbb{E}_i)$. Thus $\bigcup_{i=1}^s \mathcal{D}_{\mathcal{F}}^{\mathcal{G}}(F_i)$ is a
spanning subcomplex of $\mathcal{D}_{\mathcal{F}}^{\mathcal{G}}$.
\hfill $\square$

\begin{prop1} Let $\mathcal{F} \leq \mathcal{G}$ be finitely generated abelian
groups with $n$ generators and let $F \in K[x_1,\ldots,x_n]$ be an
$\mathcal{F}$-homogeneous polynomial. Then the simplicial complex
$\mathcal{D}_{\mathcal{F}}^{\mathcal{G}}(F)$ is a simplex.
\end{prop1}
\textbf{Proof.} The empty space is a
simplex, so it is enough to consider the case where
$\mathcal{D}_{\mathcal{F}}^{\mathcal{G}}(F)$ is not empty. Let $A,
B$ be two sets of vectors such that
$Sat(L_{\mathcal{G}})=L_{\mathbb{Z}A}$,
$Sat(L_{\mathcal{F}})=L_{\mathbb{Z}B}$ and let $T$ be the set of
vertices of $\mathcal{D}_{\mathcal{F}}^{\mathcal{G}}(F)$. Then for
every $\mathbb{E}_i \in T$ there exists a monomial $N_i={\bf
x}^{{\bf u}_i}$ in $F$ such that ${\rm deg}_{\mathbb{Z}A}({\bf
u}_i)\in relint_{\mathbb{Q}}(\sigma_{\mathcal{G}}(\mathbb{E}_i))$, see
the proof of Theorem 5.1 in \cite{KMT}. Consequently $\pi({\rm
deg}_{\mathbb{Z}A} ({\bf u}_i))$ belongs to $relint_{\mathbb{Q}}\left(
\pi(\sigma_{\mathcal{G}}({\mathbb{E}_i}))\right)$ and therefore
${\rm deg} _{\mathbb{Z}B} ({\bf u}_i) \in relint_{\mathbb{Q}}\left(
\pi(\sigma_{\mathcal{G}}({\mathbb{E}_i}))\right)$ since ${\rm deg}
_{\mathbb{Z}B} ({\bf u}_i)=\pi({\rm deg}_{\mathbb{Z}A} ({\bf
u}_i))$. But $F$ is $\mathbb{Z}B$-homogeneous, so ${\rm deg}
_{\mathbb{Z}B} ({\bf u}_i)$ is the same for all monomials in $F$.
Hence $${\rm deg}_{\mathbb{Z}B} ({\bf u}_i) \in \bigcap
_{{\mathbb{E}_i} \in T}relint_{\mathbb{Q}}\left(
\pi(\sigma_{\mathcal{G}}({\mathbb{E}_i}))\right).$$ Thus
$$\bigcap_{{\mathbb{E}_i} \in
T}relint_{\mathbb{Q}}\left(
\pi(\sigma_{\mathcal{G}}({\mathbb{E}_i}))\right) \not=\emptyset,$$
which implies that $T\in \mathcal{D}_{\mathcal{F}}^{\mathcal{G}}$
and then also $T\in \mathcal{D}_{\mathcal{F}}^{\mathcal{G}}(F)$,
since $\mathcal{D}_{\mathcal{F}}^{\mathcal{G}}(F)$ is an induced
subcomplex. Consequently
$\mathcal{D}_{\mathcal{F}}^{\mathcal{G}}(F)$ is a simplex. \hfill
$\square$\\

Combining Theorem 3.7 with Proposition 3.8 we get the following corollary:
\begin{cor1} Let $\mathcal{F} \leq \mathcal{G}$ be finitely generated abelian
groups with $n$ generators. If $F_1,\ldots,F_s$ are
$\mathcal{F}$-homogeneous polynomials and generate
$rad(I_{L_{\mathcal{G}}})$ up to radical, then $\bigcup_{i=1}^s
\mathcal{D}_{\mathcal{F}}^{\mathcal{G}}(F_i)$ is a spanning
subcomplex of $\mathcal{D}_{\mathcal{F}}^{\mathcal{G}}$ and each
$\mathcal{D}_{\mathcal{F}}^{\mathcal{G}}(F_i)$ is a simplex.
\end{cor1}
We can use Corollary 3.9 to provide a lower bound for ${\rm
ara}_{\mathbb{Z}B}(I_{L_{\mathcal{G}}})$, where $B$ is a set of vectors such that $Sat(L_{\mathcal{F}})=L_{\mathbb{Z}B}$.\\
Let $\mathcal{D}$ be a simplicial complex with vertices
$\mathcal{V}=\{v_1,\ldots,v_n \}$ and $\Omega=\{0,1,\dots
,dim(\mathcal{D})\}$. A set $\mathcal{M}=\{T_1,\dots ,T_s\}$ of
simplices of $\mathcal{D}$ is called an $\Omega$-{\em matching} in
$\mathcal{D}$ if $T_k\cap T_l=\emptyset$, for all distinct indices $k$ and $l$, see also Definition 2.1 in \cite{KT}.
 Let $supp(\mathcal{M})=\cup^{s}_{i=1}T_i$, which is a subset of
the set of vertices $\mathcal{V}$. We denote by $card(\mathcal{M})$ the
cardinality $s$ of the set $\mathcal{M}$. An $\Omega$-matching
$\mathcal{M}$ in $\mathcal{D}$ is called a
 {\em maximal} $\Omega$-matching if $supp(\mathcal{M})$ has the maximum possible
 cardinality among all $\Omega$-matchings. By $\delta(\mathcal{D})_{\Omega}$ we denote
 the minimum $card(\mathcal{M})$ among all maximal $\Omega$-matchings $\mathcal{M}$ in
 $\mathcal{D}$. For a simplicial complex $\mathcal{D}$ the number
$\delta(\mathcal{D})_{\Omega}$ is equal to
the smallest number $s$ of simplices $T_i$ of $\mathcal{D}$ such
that the subcomplex $\cup_{i=1}^sT_i$ is spanning, see Proposition 3.3 in \cite{KT}. These numbers were
introduced in \cite{KT}, where we proved that
$\delta(\mathcal{D}_{\mathcal{G}}^{\mathcal{G}})_{\Omega} \leq
{\rm ara}_{\mathcal{G}}(I_{L_{\mathcal{G}}})$. To every simplicial
complex $\mathcal{D}$ we can associate a simple graph, called the
$\{0,1\}$-{\em skeleton} of $\mathcal{D}$ and denoted by
$\mathbb{G}(\mathcal{D})$, formed by the simplices of
$\mathcal{D}$ of dimension at most $1$. The {\em complement} of
$\mathbb{G}(\mathcal{D})$, denoted by ${\overline
{\mathbb{G}(\mathcal{D} )}}$, is the graph with the same vertices
as $\mathbb{G}(\mathcal{D})$, such that there is an edge between the
vertices $v_i$ and $v_j$ if and only if there is no edge between
$v_i$ and $v_j$ in the graph $\mathbb{G}(\mathcal{D})$. Given an
integer $k$, a $k$-{\em coloring} of ${\overline
{\mathbb{G}(\mathcal{D} )}}$ is a function $c: \mathcal{V}
\rightarrow \{1,\dots ,k\}$ such that $c(v_i)\not= c(v_j)$ if the
vertices $v_i,v_j$ are joined by an edge of ${\overline
{\mathbb{G}(\mathcal{D} )}}$. The {\em chromatic number} $\gamma
({\overline {\mathbb{G}(\mathcal{D} )}})$ of ${\overline
{\mathbb{G}(\mathcal{D} )}}$ is the smallest integer $k$ such that
there is a $k$-coloring of ${\overline {\mathbb{G}(\mathcal{D}
)}}$.\\
Combining Corollary 3.9 with Corollary 2.12 in \cite{KT} we have
the following Theorem:
\begin{thm1} Let $\mathcal{F} \leq \mathcal{G}$ be finitely generated abelian
groups with $n$ generators and $B$ a set of vectors such that
$Sat(L_{\mathcal{F}})=L_{\mathbb{Z}B}$, then $$\gamma
(\overline{\mathbb{G}(\mathcal{D}_{\mathcal{F}}^{\mathcal{G}})})
\leq \delta(\mathcal{D}_{\mathcal{F}}^{\mathcal{G}})_{\Omega} \leq
{\rm ara}_{\mathbb{Z}B}(I_{L_{\mathcal{G}}}) \leq {\rm
ara}_{\mathcal{F}}(I_{L_{\mathcal{G}}}).$$
\end{thm1}
In the case that the finest $\mathcal{H}$-grading of Theorem 3.4
is given by the zero group the lower bound given by Theorem 3.10
does not provide actually any information about the arithmetical
rank of a lattice ideal. Even in this case the next theorem
provides information about the size and the complexity of the
polynomials $F_1,\ldots,F_s$ which generate
$rad(I_{L_{\mathcal{G}}})$ up to radical.
\begin{thm1} Let $I_{L_{\mathcal{G}}}$ be a lattice ideal and
$\mathcal{F}$ a specialization of $\mathcal{G}$. If
$F_{1},\ldots,F_{s}$ generate $rad(I_{L_{\mathcal{G}}})$ up to
radical, then \begin{enumerate} \item the total number of
monomials in the nonzero terms of the polynomials
$F_{1},\ldots,F_{s}$ is greater than or equal to the number of
vertices of $\mathcal{D}_{\mathcal{G}}^{\mathcal{G}}$ and \item
the total number of $\mathcal{F}$-homogeneous components in
$F_{1},\ldots,F_{s}$ is greater than or equal to
$\delta(\mathcal{D}_{\mathcal{F}}^{\mathcal{G}})_{\Omega}$.
\end{enumerate}
\end{thm1}
\textbf{Proof.} (1) Using Theorem 3.7 we take that for each vertex $\mathbb{E}_i$ of $\mathcal{D}_{\mathcal{G}}^{\mathcal{G}}$
there exists at least one monomial $N$ in a nonzero term of some $F_j$, such that $cone(N)=\sigma_{\mathcal{G}}(\mathbb{E}_i)$. The result follows.
\\ (2) Let $F_{i}({\bf b}_{i,1}),\ldots,F_{i}({\bf
b}_{i,q_{i}})$ be all the $\mathcal{F}$-homogeneous components of
$F_i$, $1 \leq i \leq s$. Then
$$rad(I_{L_{\mathcal{G}}})=rad(F_{1}({\bf
b}_{1,1}),\ldots,F_{1}({\bf b}_{1,q_{1}}),\ldots,F_{s}({\bf
b}_{s,1}),\ldots,F_{s}({\bf b}_{s,q_{s}}))$$ since
$$(F_1,\ldots,F_s) \subset (F_{1}({\bf
b}_{1,1}),\ldots,F_{1}({\bf b}_{1,q_{1}}),\ldots,F_{s}({\bf
b}_{s,1}),\ldots,F_{s}({\bf b}_{s,q_{s}})) \subset
I_{L_{\mathcal{G}}}.
$$ Thus ${\rm ara}_{\mathcal{F}}(I_{L_{\mathcal{G}}}) \leq
q_{1}+\cdots+q_{s}$ and therefore, from Theorem 3.10, we have that
$\delta(\mathcal{D}_{\mathcal{F}}^{\mathcal{G}})_{\Omega} \leq
q_{1}+\cdots+q_{s}$. \hfill $\square$\\
\newline

\section{Application}

In this section we will give an example of a toric ideal
$I_{L_{\mathbb{Z}A_{G}}}$ to explain how the techniques of the
previous sections can be applied to give lower bounds for the
$\mathcal{F}$-homogeneous arithmetical rank. For the toric ideal
$I_{L_{\mathbb{Z}A_{G}}}$ we prove that: \begin{enumerate} \item
it is not a $\mathbb{Z}B$-homogeneous, as well as
$\mathbb{Z}A_{G}$-homogeneous, set-theoretic complete
intersection, for a certain specialization $\mathbb{Z}B$ of
$\mathbb{Z}A_{G}$. \item it is not an $\mathcal{F}$-homogeneous
set-theoretic complete intersection, for infinitely many
specializations $\mathcal{F}$ of $\mathbb{Z}A_{G}$.
\end{enumerate}
One can use the techniques, based on circuits of a vector
configuration, developed in \cite{KT} to compute the simplicial
complex $\mathcal{D}_{\mathcal{G}}^{\mathcal{G}}$ and therefore
find the vertices $\{\mathbb{E}_{1},\ldots,\mathbb{E}_{f}\}$ of
the simplicial complex $\mathcal{D}_{\mathcal{F}}^{\mathcal{G}}$.
Explicitly computing the intersections of the relative interiors
of the cones $\pi(\sigma_{\mathcal{G}}({\mathbb{E}_i}))$ we obtain
the simplices of $\mathcal{D}_{\mathcal{F}}^{\mathcal{G}}$. Using
all these informations we can compute the chromatic number of the
complement of the $\{0,1\}$-skeleton of
$\mathcal{D}_{\mathcal{F}}^{\mathcal{G}}$, which provides a lower
bound for the $\mathcal{F}$-homogeneous arithmetical rank.

Let $G$ be the graph cube

$$ 
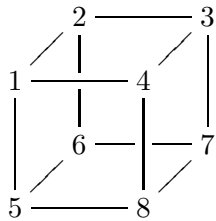
\xymatrix@!0{ & 2 \ar@{-}[rr]\ar@{-}'[d][dd] & & 3 \ar@{-}[dd]
\\
1 \ar@{-}[ur]\ar@{-}[rr]\ar@{-}[dd] & & 4 \ar@{-}[ur]\ar@{-}[dd]
\\
& 6 \ar@{-}'[r][rr] & & 7
\\
5 \ar@{-}[rr]\ar@{-}[ur] & & 8 \ar@{-}[ur] } $$ To every graph we
can assign a toric ideal in the polynomial ring with so many
variables as the edges of the graph. This toric ideal is commonly
known as the toric ideal arising from the graph $G$. More details
about toric ideals arising from finite graphs can be found in \cite{V} and in
\cite{OH}. Let $A_{G}$ be the set of all vectors ${\bf a}_{ij}={\bf
e}_i+{\bf e}_j$ such that $\{t_i,t_j\}$, $i<j$, is an edge of $G$,
where $\{{\bf e}_i \ | \ 1\leq i\leq 8\}$ is the canonical basis
of ${\mathbb{R}}^8$. Note that every vector configuration coming from a
graph is extremal. A vector configuration $A$ is called {\em
extremal} if the strongly convex rational polyhedral cone
$\sigma_{\mathbb{Z}A}$ is not generated by any proper subset of
$A$.
Let $\mathcal{F}=\mathbb{Z}B$ and $\mathcal{G}=\mathbb{Z}A_{G}$.
Consider the toric ideal
$$I_{L_{\mathcal{G}}}\subset K[x_{12}, x_{14}, x_{15}, x_{23}, x_{26},
x_{34}, x_{37}, x_{48}, x_{56}, x_{58}, x_{67}, x_{78}].$$It is
the kernel of the $K$-algebra homomorphism $$\phi: \Bbbk[x_{12},
x_{14}, x_{15}, x_{23}, x_{26}, x_{34},x_{37}, x_{48},x_{56},
x_{58}, x_{67}, x_{78}] \rightarrow \Bbbk[t_{1},\ldots,t_{8}]$$
defined by $\phi(x_{ij})={\bf t}^{{\bf a}_{ij}}$. There are $6$
circuits corresponding to the cycles of length $4$, $16$ circuits
corresponding to the cycles of length $6$ and $6$ circuits
corresponding to the cycles of length $8$. More precisely the
circuits of $A_G$ are
\begin{tabbing} $\mathcal{C}_{A_{G}}=$ \=
$\{x_{14}x_{23}-x_{12}x_{34}, x_{12}x_{56}-x_{15}x_{26},
x_{26}x_{37}-x_{23}x_{67},x_{14}x_{58}-x_{15}x_{48},$
\\ \>
$x_{37}x_{48}-x_{34}x_{78}, x_{58}x_{67}-x_{56}x_{78},
x_{23}x_{48}x_{56}-x_{26}x_{34}x_{58},$
\\ \>
$x_{14}x_{37}x_{56}-x_{15}x_{34}x_{67},x_{12}x_{37}x_{58}-x_{15}x_{23}x_{78},x_{12}x_{48}x_{67}-x_{14}x_{26}x_{78},$
\\ \>
$x_{23}x_{56}x_{78}-x_{26}x_{37}x_{58},x_{14}x_{56}x_{78}-x_{15}x_{48}x_{67},x_{26}x_{34}x_{78}-x_{23}x_{48}x_{67},$
\\ \>
$x_{15}x_{34}x_{78}-x_{14}x_{37}x_{58},x_{15}x_{26}x_{78}-x_{12}x_{58}x_{67},x_{14}x_{23}x_{78}-x_{12}x_{37}x_{48},$
\\ \>
$x_{34}x_{58}x_{67}-x_{37}x_{48}x_{56},x_{12}x_{34}x_{67}-x_{14}x_{26}x_{37},x_{15}x_{23}x_{67}-x_{12}x_{37}x_{56},$
\\ \>
$x_{12}x_{34}x_{58}-x_{15}x_{23}x_{48},x_{14}x_{26}x_{58}-x_{12}x_{48}x_{56},x_{14}x_{23}x_{56}-x_{15}x_{26}x_{34},$
\\ \>
$x_{12}x_{34}x_{56}x_{78}-x_{15}x_{23}x_{48}x_{67},x_{12}x_{34}x_{56}x_{78}-x_{14}x_{26}x_{37}x_{58},$
\\ \>
$x_{14}x_{23}x_{56}x_{78}-x_{15}x_{26}x_{37}x_{48},x_{12}x_{34}x_{58}x_{67}-x_{15}x_{26}x_{37}x_{48},$
\\ \>
$x_{14}x_{23}x_{58}x_{67}-x_{12}x_{37}x_{48}x_{56},x_{14}x_{23}x_{58}x_{67}-x_{15}x_{26}x_{34}x_{78}\}$.
\end{tabbing}
Looking at the monomials of the above circuits and considering
their minimal elements, which are the 20 monomials involved in the
first ten circuits, we get all the vertices of the simplicial
complex $\mathcal{D}_{\mathcal{G}}^{\mathcal{G}}$, see Section 4 \cite{KT}.
The complex $\mathcal{D}_{\mathcal{G}}^{\mathcal{G}}$
has $20$ vertices defined by the following sets:\\
$\mathbb{E}_1=\{14,23\}, \mathbb{E}_2=\{12,34\},
\mathbb{E}_3=\{12,56\},
\mathbb{E}_4=\{15,26\},\mathbb{E}_5=\{26,37\},$\\
$\mathbb{E}_6=\{23,67\},\mathbb{E}_7=\{14,58\},
\mathbb{E}_8=\{15,48\},\mathbb{E}_9=\{37,48\},
\mathbb{E}_{10}=\{34,78\},$\\ $\mathbb{E}_{11}=\{58,67\},
\mathbb{E}_{12}=\{56,78\},\mathbb{E}_{13}=\{23,48,56\},
\mathbb{E}_{14}=\{26,34,58\},$\\ $\mathbb{E}_{15}=\{14,37,56\},
\mathbb{E}_{16}=\{15,34,67\},\mathbb{E}_{17}=\{12,37,58\},\mathbb{E}_{18}=\{15,23,78\},$\\
$\mathbb{E}_{19}=\{12,48,67\},\mathbb{E}_{20}=\{14,26,78\}.$
\\ It
has ten 1-simplices, namely $$\{\mathbb{E}_1,\mathbb{E}_2\},
\{\mathbb{E}_3,\mathbb{E}_4\},\{\mathbb{E}_5,\mathbb{E}_6\},
\{\mathbb{E}_7,\mathbb{E}_8\}, \{\mathbb{E}_9,\mathbb{E}_{10}\},$$
$$\{\mathbb{E}_{11},\mathbb{E}_{12}\},
\{\mathbb{E}_{13},\mathbb{E}_{14}\},
\{\mathbb{E}_{15},\mathbb{E}_{16}\},
\{\mathbb{E}_{17},\mathbb{E}_{18}\},
\{\mathbb{E}_{19},\mathbb{E}_{20}\}.$$ There are no $2$-simplices.
The first ten binomials of $\mathcal{C}_{A_{G}}$ constitute a
minimal set of generators of the ideal $I_{L_{\mathcal{G}}}$ and
therefore ${\rm ara}_{\mathcal{G}}(I_{L_{\mathcal{G}}}) \leq 10$.
On the other hand the chromatic number of the complement of the
$\{0,1\}$-skeleton of $\mathcal{D}_{\mathcal{G}}^{\mathcal{G}}$ is
equal to $10$, so ${\rm ara}_{\mathcal{G}}(I_{L_{\mathcal{G}}})=
10$, see Section 4 in
\cite{KT}.\\
Consider the set of vectors
\begin{tabbing} $B=$ \= $\{{\bf b}_{12}=(5,0,3,4), {\bf b}_{14}=(3,1,5,5), {\bf b}_{15}=(4,1,4,8),
{\bf b}_{23}=(4,0,2,3),$ \\ \> ${\bf b}_{26}=(5,0,1,6), {\bf
b}_{34}=(2,1,4,4), {\bf b}_{37}=(2,1,2,6), {\bf
b}_{48}=(1,2,5,8),$
\\ \> $ {\bf b}_{56}=(4,1,2,10), {\bf b}_{58}=(2,2,4,11), {\bf
b}_{67}=(3,1,1,9), {\bf b}_{78}=(1,2,3,10)\}$.
\end{tabbing} We have that $\mathcal{F}$ is a specialization of $\mathcal{G}$, since $L_{\mathcal{G}} \subset
L_{\mathcal{F}}$. Let $\pi: \sigma_{\mathcal{G}} \rightarrow
\sigma_{\mathcal{F}}$ be the projection of cones, given by
$\pi({\bf a}_{ij})={\bf b}_{ij}$. We will compute the simplices of
$\mathcal{D}_{\mathcal{F}}^{\mathcal{G}}$. The vertices of
$\mathcal{D}_{\mathcal{F}}^{\mathcal{G}}$ are the same with
$\mathcal{D}_{\mathcal{G}}^{\mathcal{G}}$, namely
$\mathbb{E}_1,\ldots,\mathbb{E}_{20}$. There are $20$ cones of the
form $\pi(\sigma_{\mathcal{G}}({\mathbb{E}_i}))$, i.e.
$$\pi(\sigma_{\mathcal{G}}({\mathbb{E}_1}))=pos_{\mathbb{Q}}({\bf b}_{14},{\bf
b}_{23}), \ \pi(\sigma_{\mathcal{G}}({\mathbb{E}_2}))=pos_{\mathbb{Q}
}({\bf b}_{12},{\bf b}_{34}) \ \textrm{etc.}$$ By explicitly
computing the intersections of the relative interiors of the above
cones we take that the simplicial complex
$\mathcal{D}_{\mathcal{F}}^{\mathcal{G}}$ has $7$ facets:
\begin{enumerate} \item one $7$-simplex, namely $\{\mathbb{E}_{13},\mathbb{E}_{14},\mathbb{E}_{15},\mathbb{E}_{16},\mathbb{E}_{17},\mathbb{E}_{18},\mathbb{E}_{19},\mathbb{E}_{20}\}$.
\item six $1$-simplices, namely
$$\{\mathbb{E}_{1},\mathbb{E}_{2}\},
\{\mathbb{E}_{3},\mathbb{E}_{4}\},\{\mathbb{E}_{5},\mathbb{E}_{6}\},
\{\mathbb{E}_{7},\mathbb{E}_{8}\},\{\mathbb{E}_{9},\mathbb{E}_{10}\},\{\mathbb{E}_{11},\mathbb{E}_{12}\}.$$
\end{enumerate}
Note that $\Omega=\{0,1,\ldots,7\}$. We have that
$\delta(\mathcal{D}_{\mathcal{F}}^{\mathcal{G}})_{\Omega}=7$,
attained by the maximal $\Omega$-matching
$$\{\{\mathbb{E}_1,\mathbb{E}_2\},\{\mathbb{E}_3,\mathbb{E}_4\},\{\mathbb{E}_5,\mathbb{E}_6\},\{\mathbb{E}_7,\mathbb{E}_8\},\{\mathbb{E}_9,\mathbb{E}_{10}\},
\{\mathbb{E}_{11},\mathbb{E}_{12}\},
\{\mathbb{E}_{13},\mathbb{E}_{14},\ldots,\mathbb{E}_{20}\}\}.$$Remark
that $\gamma
(\overline{\mathbb{G}(\mathcal{D}_{\mathcal{F}}^{\mathcal{G}})})=7$.
Therefore $7\leq {\rm ara}_{\mathcal{F}}(I_{L_{\mathcal{G}}})$.
Moreover
\begin{tabbing} $rad(I_{L_{\mathcal{G}}})=$ \=
$rad(x_{14}x_{23}-x_{12}x_{34},x_{12}x_{56}-x_{15}x_{26},x_{26}x_{37}-x_{23}x_{67},$\\
\>
$x_{14}x_{58}-x_{15}x_{48},x_{37}x_{48}-x_{34}x_{78},x_{58}x_{67}-x_{56}x_{78},$\\
\>
$(x_{23}x_{48}x_{56}-x_{26}x_{34}x_{58})+(x_{14}x_{37}x_{56}-x_{15}x_{34}x_{67})+$\\
\>
$+(x_{12}x_{37}x_{58}-x_{15}x_{23}x_{78})+(x_{12}x_{48}x_{67}-x_{14}x_{26}x_{78}))$.
\end{tabbing} So ${\rm ara}_{\mathcal{F}}(I_{L_{\mathcal{G}}})=7$. Note also
that, since the graph $G$ is bipartite, the height of the toric
ideal $I_{L_{\mathcal{G}}}$ is equal to the number of edges minus
the number of vertices plus one, see \cite{V}, so ${\rm
ht}(I_{L_{\mathcal{G}}})=5$, which implies that
$I_{L_{\mathcal{G}}}$ is not a $\mathcal{F}$-homogeneous
set-theoretic complete intersection. Actually for any group
$\mathcal{H}$ such that $ \mathcal{F} \leq \mathcal{H} \leq
\mathcal{G}$ the toric ideal $I_{L_{\mathcal{G}}}$ is not an
$\mathcal{H}$-homogeneous set-theoretic complete intersection,
since
$$5\leq {\rm
ara}_{\mathcal{F}}(I_{L_{\mathcal{G}}})=7 \leq {\rm
ara}_{\mathcal{H}}(I_{L_{\mathcal{G}}})\leq {\rm
ara}_{\mathcal{G}}(I_{L_{\mathcal{G}}})= 10.$$ There are
infinitely many different equivalent classes of $\mathcal{H}$'s
since the rank of $\mathcal{F}$ equals $8$ and the rank of
$\mathcal{G}$ equals $5$.
\bigskip
\newline
{\bf Acknowledgment}

The authors thank the referee for his careful reading of the manuscript and his helpful remarks.

\end{document}